\documentclass{article}[12pt]
\usepackage{latexsym, amsmath, amsfonts, amssymb}
\usepackage{hyperref}
\usepackage{color}

\newcommand{\nc}{\newcommand}
\nc{\ga}{\gamma} \nc{\di}{\displaystyle}
\nc{\ek}{\protect\\[1ex]}
\nc{\N}{{\mathbb N}} \nc{\R}{{\mathbb R}} \nc{\Z}{{\mathbb Z}}
\nc{\La}{\Lambda} \nc{\la}{\lambda} \nc{\da}{\delta}
\nc{\Da}{\Delta} \nc{\na}{\nabla} \nc{\vp}{\varphi} \nc{\si}{\sigma}
\nc{\Si}{\Sigma} \nc{\al}{\alpha} \nc{\be}{\beta} \nc{\om}{\omega}
\nc{\Om}{\Omega} \nc{\pa}{\partial} \nc{\ti}{\times}
 \nc{\ve}{\varepsilon} \nc{\ra}{\rightarrow} \nc{\Ra}{\Rightarrow}
\nc{\ran}{\rangle} \nc{\lan}{\langle}

 \nc{\eq}[1]{\mbox{\rm{(\ref{E#1})}}}
\nc{\qed}{\mbox{}\nolinebreak\hfill \rule{2mm}{2mm}}
 \nc{\ha}{\frac{1}{2}}
\nc{\hra}{\hookrightarrow} \nc{\supp}{{\rm supp}\,}
\nc{\curl}{\text{curl}\,} \nc{\dense}{\hra^{\hspace{-3mm}d\,\,}}

\newtheorem{lem}{Lemma}[section]
\newtheorem{theo}[lem]{Theorem}

\renewcommand{\div}{{\rm{div}}\,}

\numberwithin{equation}{section} %\linespread{1.4}

\title{Global regularity of Leray-Hopf weak solutions to 3D Navier-Stokes
 equations}

\author{Myong-Hwan Ri\\
\small Institute of Mathematics, State Academy of Sciences, DPR
 Korea}
\date{}
  \allowdisplaybreaks[4]
\begin{document}
\bibliographystyle{alpha}
%\maketitle%

\maketitle%

\begin{abstract}
We show that any Leray-Hopf weak solution to 3D Navier-Stokes
equations  with initial values $u_0\in H^{1/2}(\R^3)$ belong to
$L^\infty(0,\infty; H^{1/2}(\R^3))$ and thus it is regular.

For the proof, first, we construct a supercritical space, the norm
of which is compared to the homogeneous Sobolev $\dot{H}^{1/2}$-norm
in that it has inverse logarithmic weight very sparsely in the
frequency domain.
 Then we obtain the energy estimates of high
frequency parts of the solution which involve the supercritical norm
on the right-hand side. Finally, we superpose the energy norm of
high frequency parts of the solution to get estimates of the
critical norms of the weak solution via the re-scaling argument.

\end{abstract}

%%%%%%%%%%%%%%%%%%%%%%%%%%%%%%%%%%%%%%%%%%%%%%%%%%%%

\let\thefootnote\relax\footnote{\hspace{-0.3cm}
E-mail address\,$:$ math.inst@star-co.net.kp (Myong-Hwan Ri)}

%%%%%%%%%%%%%%%%%%%%%%%%%%%%%%%%%%%%%%%%%%%%%%%%%%%%
%%%%%%%%%%%%%%%%%%%     1 Introduction   %%%%%%%%%%%%%%%%%%%%%%%

\section{Introduction and main result}
Let us consider the Cauchy problem for the 3D Navier-Stokes
equations:
\begin{equation}
\label{E1.1}
\begin{array}{rl}
     u_t-\nu\Da u  + (u\cdot\na)u+\na p  = 0 \,\, &\text{in }(0,\infty)\ti\R^3,\ek
      \div u = 0 \,\, &\text{in }(0,\infty)\ti\R^3,\ek
u(0,x)=u_0\, &\text{in }\R^3,
\end{array}
\end{equation}
where $u$, $p$ are velocity and pressure, respectively, and $\nu$ is
the kinematic viscosity of the observed viscous incompressible
fluid.

Since Leray \cite{Le34} proved existence of a global weak solution
(Leray-Hopf weak solution) to \eq{1.1} such that the weak solution
$$
\label{E1.3} u\in L^2(0,\infty; H^1(\R^3))\cap L^\infty(0,\infty; L^2(\R^3))
 $$
 satisfies \eq{1.1} in a weak sense and the {\it energy inequality}
\begin{equation}
\label{EEI}
\begin{array}{l}\di\frac{1}{2}\|u(t)\|_2^2+\nu\int_0^t\|\na
u(\tau)\|_2^2\,d\tau\leq \frac{1}{2}\|u_0\|_2^2,\; \forall t\in
(0,\infty),
 \end{array}
\end{equation}
the problem of global regularity of the weak solution has been a
consistent key issue.
% A weak solution $u$  to \eq{1.1} is called regular (or
%strong) in  $(0,T]$, $0<T<\infty$, if
%\begin{equation}
%\label{E1.4n}
% u\in L^\infty(0,T; H^1(\R^3))\cap L^2(0,T; H^2(\R^3)).
%  \end{equation}

There is a great number of articles on scaling-invariant regularity
criteria for weak solutions to \eq{1.1}.
 We recall that the Navier-Stokes equations are
invariant by the scaling $u_\la(t,x)\equiv \la u(\la^2t,\la
x),\la>0,$ and
 a critical space for the equations is the space whose norm is invariant with respect to
the scaling
 $$ u(x)\mapsto\la u(\la x), \; \la>0.$$
The following embedding holds
 between critical spaces for the 3D Navier-Stokes equations:
$$\begin{array}{l} \dot{H}^{1/2}\hra   L^3\hra
\dot{B}^{-1+3/p}_{p,\infty}\hra
  BMO^{-1}
  \hra   \dot{B}^{-1}_{\infty,\infty},\; 3 \leq p<\infty.
  \end{array}$$
 Eskauraiza, Seregin and Sv\'erak in \cite{ESS03}(2003) proved that a Leray-Hopf weak solution
 $u$ to \eq{1.1} is regular in $(0,T]$ if
$$u\in L^\infty(0,T; L^3(\R^3)),
 $$
  employing the backward uniqueness property of parabolic equations.
By developing a profile decomposition technique and using the method of ``critical
elements" developed in \cite{KeMe06}-\cite{KeMe10}, Gallagher, Koch and Planchon proved
for  a mild solution $u$ to \eq{1.1} with initial values in $L^3(\R^3)$ that a potential
 singularity at $t=T$ implies $\lim_{t\ra
 T-0}\|u(t)\|_3=\infty$ in \cite{GKB13}(2013), and extended the result
  from $L^3(\R^3)$ to wider critical Besov spaces
$\dot{B}^{-1+3/p}_{p,q}(\R^3)$, $3<p,q<\infty$
   in \cite{GKB16}(2016).
  In \cite{Ri23}(2023), the author proved that the above results can be extended
 to the largest critical space $\dot{B}^{-1}_{\infty,\infty}(\R^3)$, that is,
 a Leray-Hopf weak solution to \eq{1.1} satisfying
$$u\in L^\infty(0,T; \dot{B}^{-1}_{\infty,\infty}(\R^3))
 $$
 is regular in $(0,T]$ by developing a method of energy superposition of high frequency parts
 of the weak solution.

Concerning supercritical regularity criteria, we mention recent
results \cite{BaPr21} and \cite{BaPr22} by Barker and Prange that a
Leray-Hopf weak solution $u$ satisfying the local energy inequality
is regular in $(0,T]$ if
$$\sup_{0<t<T}\int_{\R^3} \frac{|u(x,t)|^3}
{\Big(\log\log\log\big(\big(\log(e^{e^{3e^e}}+|u(x,t)|)\big)^{1/3}\big)\Big)^\theta}dx<\infty,\forall
\theta\in (0,1),$$ based on a Tao's result \cite{Tao19, Tao21} on
quantitative bound of the critical $L^3$-norm near a possible
blow-up epoch for a weak solution to \eq{1.1}; we also mention Pan
\cite{Pan16} and Seregin \cite{Se22} where logarithmically
supercritical regularity criteria for axisymmetric suitable weak
solutions to \eq{1.1} are obtained.

\par\medskip
In this paper, we prove that a Leray-Hopf weak solution to \eq{1.1}
with initial value in the critical space $\dot{H}^{1/2}(\R^3)$ is
regular.

\par\medskip
The main result of the paper is as follows.
\begin{theo}
\label{T1.1} {\rm Let $u$ be a Leray-Hopf weak solution to \eq{1.1}
with $u_0\in H^{1/2}(\R^3)$, $\div u_0=0$.
 Then,
$$u\in
L^\infty(0,\infty;H^{1/2}(\R^3))\cap L^2(0,\infty;H^{3/2}(\R^3))
$$
and the estimate
 \begin{equation}
 \label{E1.2}
 \di\|u(t)\|_{H^{1/2}}^2+\nu
\int_0^{t} \|u(\tau)\|_{H^{3/2}}^2\,d\tau
  \leq  2\|u_0\|_{H^{1/2}}^2,\quad \forall t\in (0,\infty),
  \end{equation}
  holds. Thus $u$ is globally regular.
 }
\end{theo}

\par\medskip
We use the following notations. The sets of all natural numbers and
all integers are denoted by $\N$ and $\Z$, respectively. The usual
$L^p$-norm for $1\leq p\leq \infty$ is denoted by $\|\cdot\|_p$.
Three-dimensional Fourier transform of a function $u$ is given by
${\cal F}u\equiv \hat{u}:=(2\pi)^{-3/2}\int_{\R^3}{u}(t,
x)e^{ix\cdot\xi}\,dx$ and ${\cal F}^{-1}u\equiv \check{u}$. For
$s\in\R$, $\dot{H}^s(\R^3)$ and $H^s(\R^3)$ stand for the
homogeneous and inhomogeneous Sobolev spaces with norms
$$\|u\|_{\dot{H}^s(\R^3)}\equiv
 \| |\xi|^s\hat{u}(\xi)\|_2\quad\text{and}\quad \|u\|_{H^s(\R^3)}\equiv
 \|(1+|\xi|^2)^{s/2}\hat{u}(\xi)\|_2,$$ respectively.
  We use the notation
\begin{equation} \label{E2.1}
\begin{array}{l}
 u^k(t,x):=\frac{1}{(2\pi)^{3/2}}\int_{|\xi|\geq k}{\mathcal F}{u}(t,
\xi)e^{ix\cdot\xi}\,d\xi,\ek u_k:=u-u^k,\;
  u_{h,k}:=u^{h}-u^{k}\quad \text{for } 0\leq h<k.
  \end{array}
\end{equation}
We do not distinguish between the spaces of vectorial functions and
scalar functions.

Let us explain the ideas for the proof of Theorem \ref{T1.1}. As is
well-known, the cancelation property $((v\cdot\na)u,u)=0$ for
suitably smooth $u,v$ with $\div v=0$ is the key to obtain global
existence of the Leray-Hopf weak solutions to \eq{1.1}. However,
this innermost property was not used, so long as the author knows,
in most previous works for regularity. Indeed, as long as the
convection term $(u\cdot\na)u$ is treated as merely a quadratic
nonlinear term, whatever improvement of its estimate could be made,
regularity of the weak solution is obtained under a smallness
condition of a critical norm of initial values or under additional
scaling-invariant conditions on the weak solution itself.

 In order to circumvent such situations, we use essentially the
 cancelation property by testing the momentum equation with high frequency parts $u^k$,
 $k\in\N$. In that case, by the cancelation $((u\cdot\na)u^k,u^k)=0$,
we are led to the kind of estimates
$$\frac{d}{2dt}\|u^k\|_2^2+\nu\|\na u^k\|_2^2\leq c(k)\|u\|_{X_1}\|\na
  u^{k/2}\|_2^2,\;\forall k\in\N,$$ with a supercritical space $X_1$,
  where $c(k)$ depends on the topology of $X_1$.
Though this estimate, at the first glance, still seems of no special
use due to the factor $c(k)$, we pay attention to the fact that
$$\sum_{k\in\N}\|u^k\|_2^2\sim
\|u^1\|_{\dot{H}^{1/2}}^2,\quad \sum_{k\in\N}\|\na u^k\|_2^2\sim
\|\na u^1\|_{\dot{H}^{1/2}}^2$$ and
$$\sum_{k\in\N}c(k)\|\na  u^{k/2}\|_2^2 \lesssim \|\na u\|_2^2
+\sum_{n\in\N}\Big(\sum_{k=1}^n c(k)\Big)\|\na u_{n,n+1}\|_2^2.$$
Thus, if we can construct a supercritical space $X_1$ such that the
averaging condition $\sum_{k=1}^n c(k)\lesssim n$ is satisfied, then
we are led to
$$\frac{d}{2dt}\|u^1\|_{\dot{H}^{1/2}}^2+\nu\|\na u^1\|_{\dot{H}^{1/2}}^2\lesssim
 \|u(t)\|_{X_1}(\|\na u\|_2^2+ \|\na u^{1}\|_{\dot{H}^{1/2}}^2),$$
which can finally yields the estimate of
$\|u(t)\|_{\dot{H}^{1/2}}^2$ provided uniform smallness of
$\|u(\tau)\|_{X_1}$ in $\tau \in [0,t]$ is guaranteed by re-scaling.

We successfully construct a supercritical space $X_1$ with the
above-mentioned averaging condition and the uniform smallness
property by re-scaling, which has a very sparse inverse logarithmic
weight in the frequency domain compared to the critical
$\dot{H}^{1/2}$-norm, Section 2.

The proof of Theorem \ref{T1.1} is given in Section 3.

\section{A frequency-weighted scaling-variant space}

In this section, we introduce a supercritical space, the norm of
which is weaker than that of $\dot{H}^{1/2}(\R^3)$, being very close
to it.
\par\medskip
Let
$$\Da_j u= (\chi_j(\xi)\hat{u}(\xi))^\vee\quad\text{for }j\in\Z, u\in {\cal
S}'(\R^3) \;\text{with }\hat{u}\in L_{loc}^1(\R_\xi^3),
$$
where $\chi_j$ is the characteristic function of the set
$\{\xi\in\R^3:  2^{j-1}\leq |\xi|< 2^{j}\}$.
 For $s\in\R$, by Plancherel's theorem,
$$\begin{array}{l}
 2^{2s(j-1)}\|\Da_j u\|_2^2 = 2^{2s(j-1)}\int_{2^{j-1}\leq |\xi|<2^j}|\hat{u}(\xi)|^2\,d\xi\ek
 \qquad \leq \int_{2^{j-1}\leq |\xi|<2^j} |\xi|^{2s}|\hat{u}(\xi)|^2\,d\xi\leq 2^{2sj}
    \|\Da_ju\|_2^2, \forall u\in \dot{H}^s(\R^3),
\end{array}$$ and the homogeneous Sobolev norm $\|u\|_{\dot{H}^s(\R^3)}$ is
equivalent to the norm
\begin{equation}
 \label{E3.1}
\big(\sum_{j\in\Z}\|2^{sj}\Da_j u\|^2_{L^2(\R^3)}\big)^{1/2}.
 \end{equation}
\medskip\indent
 Let an infinite sequence
$\{a(j)\}_{j\in\Z}\subset \R$ be such that
 \begin{equation}
 \label{E2.2}
a(j):=\left\{\begin{array}{cl}  \log_2 j & \text{if
}j=i+2^{2^k}\quad\text{for some }k\in\N, i\in\Z, -k\leq i\leq k,\ek
                1      & \text{else}.
             \end{array}
      \right.
      \end{equation}
Let us define the space $X_1$ by
 \begin{equation}
 \label{E2.3n}
 \begin{array}{l}
X_1:=\{v\in \mathcal S'(\R^3): \hat{v}(\xi)\in
L^1_{loc}(\R_\xi^3),\;\{ 2^{j/2}a^{-1}(j)\|\Da_j v\|_2\}_{j\in\Z}\in
l_2\},\ek
 \|v\|_{X_1}:=\Big\|\{2^{j/2}a^{-1}(j)\|\Da_j v\|_2\}_{j\in\Z}\Big\|_{l_2},
 \end{array}
 \end{equation}
 where the sequence $\{a(j)\}$ is given by \eq{2.2}.

Obviously, it holds that $a(j)=1$ for $j\leq 0$ and
$$\dot{H}^{1/2}(\R^3)\hra X_1,\quad \|v\|_{X_1}\leq \|v\|_{\dot{H}^{1/2}(\R^3)},\forall v\in \dot{H}^{1/2}(\R^3).$$

We pursue supercritical properties of the space $X_1$ via the next
lemma.

\begin{lem}
 \label{L4.2n}
  {\rm
 (i) Let $v\in \dot{H}^{1/2}(\R^3)$. Then,
$$\begin{array}{l}\forall \ve>0, \exists l_0:=\max\{M+1,
\log_2\big[\frac{\log_2M\|v\|_{X_1}}{\ve}\big]+1\}>0,\ek
  \|\la v(\la\cdot)\|_{X_1}\leq \ve,\quad  \forall
\la=2^{2^{2^l}} (l\geq l_0),
\end{array}$$
 where $M>0$ is such that
\begin{equation}
 \label{E4.4}
\sum_{|j|\geq M} 2^j \|\Da_j v\|_2^2\leq \frac{\ve^2}{2}.
 \end{equation}

(ii) Let $u\in C\big([0,T],\dot{H}^{1/2}(\R^3)\big)$ with
$0<T<\infty$. Then, for any sufficiently small $\ve>0$ there exists
$l_0>0$ such that
$$\begin{array}{l}
  \|\la u(t, \la\cdot)\|_{X_1}\leq \ve,\quad  \forall t\in [0,T],\;\forall
\la=2^{2^{2^l}} (l\geq l_0).
\end{array}$$

}
\end{lem}
 %%%
 %%%
 {\bf Proof:} -- {\it Proof of (i)}:
Let $v\in \dot{H}^{1/2}(\R^3)$ and fix $\ve>0$ arbitrarily. Then, in
view of
$$ \|v\|_{\dot{H}^{1/2}}\sim \Big\|\{2^{j/2}\|\Da_j v\|_2\}_{j\in\Z}\Big\|_{l_2},$$
  there is $M=M(\ve,v)\in \N$ such that
 \begin{equation}
 \label{E4.4}
\sum_{|j|\geq M} 2^j \|\Da_j v\|_2^2\leq \frac{\ve^2}{2}.
 \end{equation}

 In the proof, for the moment, we use a short notation $v_\la:=\la v(\la\cdot)$.
  By definition of the space $X_1$ we have
 $$\|v_\la\|_{X_1}= \Big\|\{2^{j/2}a^{-1}(j)\|\Da_j v_\la\|_2\}_{j\in \Z}\Big\|_{l_2}.$$
  Observe that $({\mathcal F} v_\la)(\xi)= \la ({\mathcal F}
v(\la\cdot))(\xi)=\la^{-2}(\mathcal Fv)(\frac{\xi}{\la})$ for all
$\la>0$. Then, for all $j\in\Z$
  $$\begin{array}{l}\Da_j v_\la (x)
   =\di\frac{1}{(2\pi)^{3/2}}\int_{2^{j-1}\leq |\xi|<2^j} {\cal F}v_\la(\xi)
e^{ix\cdot\xi}\,d\xi\\[3ex]
  \qquad =\di\frac{\la^{-2}}{(2\pi)^{3/2}}\int_{2^{j-1}\leq |\xi|<2^j} {\cal
F}v(\frac{\xi}{\la}) e^{ix\cdot\xi}\,d\xi\\[3ex]
 \qquad =\di\frac{\la}{(2\pi)^{3/2}}\int_{\la^{-1}2^{j-1}\leq |\eta|<\la^{-1}2^j} {\cal
F}v(\eta) e^{i \la x\cdot\eta} \,d\eta.
  \end{array}$$
Hence, if $\la\equiv 2^{2^{2^l}}$,  $l\in\N$, then
 \begin{equation}
 \label{E2.6n}
 \begin{array}{l}\Da_j v_{\la}(x)=\di\frac{\la}{(2\pi)^{3/2}}
  \int_{2^{j-2^{2^l}-1}\leq |\eta|<2^{j-2^{2^l}}} {\cal
F}v(\eta) e^{i \la x\cdot\eta}\,d\eta \\[3ex]
 \qquad =\la (\Da_{j-2^{2^l}}v)(\la x),
 \end{array}
 \end{equation}
 yielding
 $$\begin{array}{l}2^{j/2}\Da_j v_{\la}(x)=2^{j/2+2^{2^l}} (\Da_{j-2^{2^l}}v)(\la x),\;\forall j\in\Z.
 \end{array}$$
Note that
 $$\|(\Da_{j-2^{2^l}}v)(\la \cdot)\|_2=\la^{-3/2}\|\Da_{j-2^{2^l}}v\|_2.$$
Therefore, by \eq{2.6n} we have
% \begin{equation}
% \label{E2.7}
%2^{-j}a^{-1}(j)\|\Da_j
%v_{\la_0}\|_\infty=\frac{a(j-2^{2^l})}{a(j)}2^{-j+2^{2^l}}a^{-1}(j-2^{2^l})
% \|\Da_{j-2^{2^l}}v\|_\infty,\;\forall j\in\Z,
% \end{equation}
 \begin{equation}
 \begin{array}{rcl}
\di 2^{j/2}\|\Da_j v_{\la}\|_2 &= & 2^{j/2+2^{2^l}}
\|(\Da_{j-2^{2^l}}v)(\la \cdot)\|_2\ek
 &= &\di 2^{\frac{1}{2}(j-2^{2^l})}\|\Da_{j-2^{2^l}}v\|_2,\;\forall j\in\Z,
 \end{array}
 \end{equation}
and, in view of \eq{4.4},
$$\begin{array}{l}
 \|v_{\la}\|_{X_1}^2=\di\Big\|\big\{2^{j/2}a^{-1}(j)\|\Da_j v_{\la}\|_2\big\}_{j\in\Z}\Big\|_{l_2}^2\\[2ex]
 =\di\Big\|\big\{2^{\frac{1}{2}(j-2^{2^l})}a^{-1}(j)
 \|\Da_{j-2^{2^l}}v\|_2\big\}_{j\in\Z}\Big\|_{l_2}^2\\[2ex]
 =\di\Big\|\big\{\frac{1}{a(j+2^{2^l})}2^{j/2}
 \|\Da_{j}v\|_2\big\}_{j\in\Z}\Big\|_{l_2}^2\\[2ex]
 \leq \di\Big\|\big\{\frac{1}{a(j+2^{2^l})}2^{j/2}
 \|\Da_{j}v\|_2\big\}_{|j|\leq M}\Big\|_{l_2}^2+\Big\|\big\{2^{j/2}
 \|\Da_{j}v\|_2\big\}_{|j|> M}\Big\|_{l_2}^2\\[2ex]
  \di\leq \frac{\ve^2}{2}+\di\Big\|\big\{\frac{a(j)}{a(j+2^{2^l})}2^{j/2}a^{-1}(j)
 \|\Da_{j}v\|_2\big\}_{-M\leq j\leq M}\Big\|_{l_2}^2.
 \end{array}$$
 Here, for $l\geq M+1$, by definition of $a(j)$ we have
 $$\begin{array}{l}
  \di\big\|\{\frac{a(j)}{a(j+2^{2^l})}2^{j/2}a^{-1}(j)
 \|\Da_{j}v\|_2\}_{-M\leq j\leq M}\big\|_{l_2}\ek\qquad
=  \di\big\|\{\frac{a(j)}{\log_2(j+2^{2^l})}2^{j/2}a^{-1}(j)
 \|\Da_{j}v\|_2\}_{-M\leq j\leq M}\big\|_{l_2}\ek\qquad
 \leq \di\sup_{-M\leq j\leq M}\frac{\log_2M}{\log_2(j+2^{2^l})}\cdot
 \big\|\{2^{j/2}a^{-1}(j)\|\Da_{j}v\|_2\}_{-M\leq j\leq M}\big\|_{l_2}
 \\[3ex]
 \qquad \leq \di\frac{\log_2M}{2^{l-1}}\|v\|_{X_1}.
  \end{array}$$
Therefore, for all
 \begin{equation}
 \label{E2.11}
l\geq l_0:=\max\{M+1,
\log_2\big[\frac{\log_2M\|v\|_{X_1}}{\ve}\big]+1\}, \quad
\la=2^{2^{2^l}},
 \end{equation}
  we have
$$ \|v_{\la}\|_{X_1}^2\leq \frac{\ve^2}{2}+\di\Big(\frac{\log_2M}{2^{l-1}}\|v\|_{X_1}\Big)^2
  \leq \ve^2.
 $$

%Thus, in view of arbitrariness of $\ve>0$, the lemma is proved.
\par\bigskip
{\it -- Proof of (ii)}: Let $u\in C([0,T];\dot{H}^{1/2}(\R^3))$ and
put
$$f(t,s):=\|u_{1/s}(t)\|_{\dot{H}^{1/2}}+\|u^s(t)\|_{\dot{H}^{1/2}},\quad (t,s)\in [0,T]\ti [2,\infty).$$
Then,  $u\in C([0,T];\dot{H}^{1/2}(\R^3))$ implies that $f\in
BC([0,T]\ti [2,\infty))$ and for all $t\in [0,T]$ the mapping $s \to
f(t,s)$ is nonincreasing in $s$ and $f(t,s)\rightarrow 0$ as $s\ra
\infty$. Therefore, for any sufficiently small $\ve >0$ and $t\in
[0,T]$ there exists $s=s(t)$ satisfying $f(t,s(t))=\ve/2$, where
 one may choose such $s(t)$ that
 \begin{equation}\label{E2.12}
\exists N>0,\; \forall t\in [0,T]:  \quad |s(t)|\leq N
 \end{equation}
  holds. In
fact, if \eq{2.12} does not hold true,  there is a sequence
$\{t_n\}\subset [0,T]$ satisfying $f(t_n, s(t_n))=\ve/2$ and
$s(t_n)\ra\infty$ as $n\ra\infty$. Without loss of generality, we
may assume $t_n\ra \tilde{t}$ for some $\tilde{t}\in [0,T]$ for
$n\ra \infty$. Then, in view of $u\in C([0,T],\dot{H}^{1/2}(\R^3))$,
we are led to
$$\begin{array}{rl}
 \di\frac{\ve}{2}&=f(t_n, s(t_n))\leq |f(t_n, s(t_n))-f(\tilde{t}, s(t_n))|+f(\tilde{t},
 s(t_n))\ek
 &\leq \big|\|[u(t_n)]_{1/s(t_n)}\|_{{\dot{H}^{1/2}}}-\|[u(\tilde{t})]_{1/s(t_n)}\|_{\dot{H}^{1/2}}\big|\ek
&\qquad
+\big|\|[u(t_n)]^{s(t_n)}\|_{{\dot{H}^{1/2}}}-\|[u(\tilde{t})]^{s(t_n)}\|_{\dot{H}^{1/2}}\big|+f(\tilde{t},
 s(t_n))\ek
 &\leq \|[u(t_n)-u(\tilde{t})]_{1/s(t_n)}\|_{\dot{H}^{1/2}}+\|[u(t_n)-u(\tilde{t})]^{s(t_n)}\|_{\dot{H}^{1/2}}+f(\tilde{t},
 s(t_n))\ek
 &\leq \sqrt{2}\|u(t_n)-u(\tilde{t})\|_{\dot{H}^{1/2}}+f(\tilde{t},
 s(t_n))\ek
 &\ra 0 \quad (n\ra\infty),
  \end{array}$$
which is a contradiction.

Consequently, if we put
$$M:=[\log_2N]+2$$
with $N$ in \eq{2.12}, then, in view of $N\leq 2^{M-1}$, we have
$$\begin{array}{rcl}
 \di\sum_{|j|\geq M} 2^j \|\Da_j u(t)\|_2^2 &= & 2\di\sum_{|j|\geq M} 2^{j-1} \|\Da_j
 u(t)\|_2^2\ek
&\leq &
2\big(\|[u(t)]_{2^{-M+1}}\|_{\dot{H}^{1/2}}^2+\|[u(t)]^{2^{M-1}}\|_{\dot{H}^{1/2}}^2\big)\ek
&\leq &
2\big(\|[u(t)]_{1/N}\|_{\dot{H}^{1/2}}^2+\|[u(t)]^{N}\|_{\dot{H}^{1/2}}^2\big)\ek
&\leq & \di\frac{\ve^2}{2}
 \end{array}$$
and thus \eq{4.4} holds true for $u(t)$ uniformly for all $ t\in
[0,T]$.

Thus,  for
$$l\geq l_0:=\max\{M+1,
\log_2\big[\frac{\log_2M\|u\|_{C([0,T],
\dot{H}^{1/2})}}{\ve}\big]+1\},$$ see \eq{2.11},  we have
$$ \|(u(t))_{\la}\|_{X_1}^2\leq \frac{\ve^2}{2}+\di\Big(\frac{\log_2M}{2^{l-1}}\|u(t)\|_{X_1}\Big)^2
  \leq \ve^2,\quad \forall t\in [0,T],
 $$
by the above proved fact (i) and $\|u(t)\|_{X_1}\leq
\|u(t)\|_{\dot{H}^{1/2}}$.

 \hfill\qed

 %%%%%%
%%%%
 %%%%%%%%
\section{Proof of the main theorem}

\begin{lem}
 \label{L4.1}
 {\rm
 Let
  \begin{equation}
  \label{E3.2}
 b(j) \equiv 2^{-j-1}\sum_{i=1}^{j}
 2^{i}a(i),\;j\in\N.
  \end{equation}
   Then, for all $j\in\N$
$$b(j)\leq \left\{\begin{array}{cl}  \log_2 j
     & \text{if }-k+2^{2^k}\leq j\leq 2^{2^k}+2k\quad\text{for some }k\in\N,\ek
                2      & \text{else}.
             \end{array}
      \right.
      $$
      }
\end{lem}
 {\bf Proof:}
 Note that $a(j)\geq 1$,  and
$b(j+1)=\frac{1}{2}(b(j)+a(j+1))$ for all $j\in\N$.

 First we prove the lemma for $j=1\sim  2^{2^1}=4$. We have
 \begin{equation}
 \label{E2.13}
\begin{array}{l}
  b(1)=2^{-2}\cdot 2\cdot a(1)=\frac{1}{2}<a(1)=1,\ek
  b(2)=2^{-3}\cdot (2+2^2)=\frac{3}{4}<a(2),\ek
  b(3)=\frac{b(2)+a(3)}{2}=\frac{3+4\log_2 3}{8}\leq \log_2 3=a(3),\ek
b(4)=\frac{b(3)+a(4)}{2}\leq 2= a(4).
 \end{array}
 \end{equation}

Next, let us fix any $k\in\N$ and prove the lemma for $2^{2^k}<
j\leq 2^{2^{k+1}}$. Let $j_1=2^{2^k}$ and assume in view of
\eq{2.13} that $b(j_1)\leq a(j_1)$ which, at the moment, is
satisfied for $k=1$.
 Then, for
all $l=1,2,\ldots,k$ we have
$$\begin{array}{l}
  b(j_1+1)=\frac{b(j_1)}{2}+\frac{a(j_1+1)}{2}\leq \frac{a(j_1)}{2}+\frac{a(j_1+1)}{2}\leq a(j_1+1)=\log_2(j_1+1),\ek
   b(j_1+2)=\frac{b(j_1+1)}{2}+\frac{\cdot a(j_1+2)}{2}\leq \frac{a(j+1)}{2}+\frac{a(j_1+2)}{2}
    \leq a(j_1+2)=\log_2(j_1+2),\ek
     \hspace{3cm}\cdots\quad\cdots\quad\cdots \ek
  b(j_1+l)=\frac{b(j_1+l-1)}{2}+\frac{a(j_1+l)}{2}\leq a(j_1+l)=\log_2(j_1+l).
  \end{array}$$
Moreover, for $l\in\N $ with $1\leq l \leq 2^{2^{k+1}}-(2k+2^{2^k})$
we have
$$\begin{array}{l}
  b(j_1+k+1)=\frac{b(k+2^{2^k})}{2}+\frac{a(k+2^{2^k}+1)}{2}\leq \frac{a(k+2^{2^k})}{2}+\frac{1}{2}
   \leq \frac{\log_2(k+2^{2^k})}{2}+\frac{1}{2},\ek
  b(j_1+k+2)=\frac{b(k+2^{2^k}+1)}{2}+\frac{a(k+2^{2^k}+2)}{2}\leq \frac{\log_2(k+2^{2^k})}{4}
       +\frac{1}{4}+\frac{1}{2},\ek
     \hspace{3cm}\cdots\quad\cdots\quad\cdots \ek
  b(j_1+k+l)=\frac{b(k+2^{2^k}+l-1)}{2}+\frac{a(k+2^{2^k}+l)}{2}\leq
      \frac{\log_2(k+2^{2^k})}{2^{l}}+\frac{1}{2^{l}}+\cdots\frac{1}{4}+\frac{1}{2}.
  \end{array}$$
Hence, if $1\leq l <\log_2\log_2(k+2^{2^k})$, then we have
 $$\begin{array}{l} b(j_1+k+l)\leq \log_2(k+2^{2^k})
  \end{array}$$
  and, if $\log_2\log_2(k+2^{2^k})\leq l \leq 2^{2^{k+1}}-(2k+2^{2^k})$, then we have
 $$\begin{array}{l} b(j_1+k+l)\leq 2.
  \end{array}$$
Note that $k<\log_2\log_2(k+2^{2^k})<k+1$ for any $k\in\N$.
Consequently,
 \begin{equation}\label{E2.10}
\begin{array}{l}
 b(j_1+k+l)\leq \log_2(k+2^{2^k}), \quad \text{for }l=1,\ldots, k,\ek
b(j_1+k+l)\leq 2, \quad\text{for }l=k+1,\ldots,
2^{2^{k+1}}-(2k+2^{2^k}).
  \end{array}
  \end{equation}

Finally, we shall prove
 \begin{equation}
 \label{E2.14}
b(j)\leq \log_2 j\quad\text{for }2^{2^{k+1}}-k+1\leq j\leq
2^{2^{k+1}}.
 \end{equation}
In fact, we have
$$\begin{array}{l}
b(2^{2^{k+1}}-k+1)=\frac{b(2^{2^{k+1}}-k)}{2}+\frac{a(2^{2^{k+1}}-k+1)}{2}\ek
  \qquad\qquad \leq 1+\frac{\log_2(2^{2^{k+1}}-k+1)}{2}\leq \log_2(2^{2^{k+1}}-k+1), \ek
 b(2^{2^{k+1}}-k+2)=\frac{b(2^{2^{k+1}}-k+1)}{2}+\frac{a(2^{2^{k+1}}-k+2)}{2}\ek
 \qquad \qquad \leq
  \frac{\log_2(2^{2^{k+1}}-k+1)}{2}+\frac{\log_2(2^{2^{k+1}}-k+2)}{2}
   \leq \log_2(2^{2^{k+1}}-k+2),\ek
 \hspace{3cm}\cdots\quad\cdots\quad\cdots \ek
 b(2^{2^{k+1}})=\frac{b(2^{2^{k+1}}-1)}{2}+\frac{a(2^{2^{k+1}})}{2}\ek
 \qquad \qquad \leq
  \frac{\log_2(2^{2^{k+1}})}{2}+\frac{\log_2(2^{2^{k+1}}-1)}{2}
   \leq 2^{k+1}=\log_2{2^{2^{k+1}}}=a(2^{2^{k+1}}).
\end{array}$$
Here, in particular, $b(2^{2^{k+1}})\leq a(2^{2^{k+1}})$ together
with $b(4)\leq a(4)$  implies by induction that the assumption
$b(j_1)\leq a(j_1)$ is satisfied for all $k\in\N$.

Thus, by \eq{2.13} $\sim$ \eq{2.14}, the proof of the lemma is
complete.

 \hfill\qed

In the below we use the notation
 \begin{equation}
 \label{E4.1}
 S(n)=\{l\in\N:\; l\leq n, -k+2^{2^k}\leq l\leq 2^{2^k}+2k,\exists k\in\N\},\;n\in\N.
  \end{equation}
  %%%
  %%%
\begin{lem}
 \label{L4.2}
 {\rm
The cardinal number $|S(n)|$ of the set $S(n)$ is estimated by
 $$|S(n)|\leq \di\frac{(3\log_2\log_2n +5)\log_2\log_2n}{2} .$$
 }
\end{lem}
 {\bf Proof:} Suppose that  $2^{2^{k-1}}\leq n <2^{2^k}$ for some $k\in\N$.
  Then, by definition of $S(n)$, we have
$$\begin{array}{rcl}
 |S(n)| & \leq  & \di\sum_{l=1}^{k-1} (3l+1)=\frac{(3k+2)(k-1)}{2}
 \\[3ex]
 &< & \di\frac{(3\log_2\log_2n +5)\log_2\log_2n}{2} .
\end{array}
$$
Thus, the lemma is proved. \hfill\qed

\par\bigskip\noindent
 {\bf Proof of Theorem \ref{T1.1}:} The proof is divided in two
 steps.

{\it Step 1. Energy estimate of high frequency parts:}

Let  $(0,T)$ be the maximal interval where a Leray-Hopf weak
solution $u$ to \eq{1.1} is regular.

Recall the notation \eq{2.1}. By $L^2$-scalar product to \eq{1.1}
with $u^k(t)$, $k\geq 1$, $t\in (0,T)$, we have
 \begin{equation}
 \label{E2.18}
\begin{array}{l}
\frac{d}{2dt}\|u^k\|_2^2+\nu\|\na u^k\|_2^2=-((u\cdot\na)u_k,
u^k)\ek
 \quad =-((u_k\cdot\na)u_k+(u^k\cdot\na)u_k, u^k)\ek
 \quad
=-((u_k\cdot\na)u_{k/2,k},u^{k})-((u_{k/2,k}\cdot\na)u_{k/2},u^{k})-((u^k\cdot\na)u_k,u^k)
 \;\text{in } (0,T)
\end{array}
 \end{equation}
  in view of  $u_k=u_{k/2}+u_{k/2,k}$ and
$$ \begin{array}{l}
 (u_k, u^l)_2=0,\quad \supp{\widehat{(u_ku_l)}}\subset \{\xi\in\R^3: |\xi|\leq k+l\},\;\quad k\leq l<\infty.
 \end{array}
 $$

 In the right-hand side of \eq{2.18} we have
$$\begin{array}{l}
|((u_{k/2,k}\cdot\na)u_{k/2},u^{k})+(u^k\cdot\na)u_k,u^k)|
 \leq  (\|\na u_{k/2}\|_\infty+\|\na u_k\|_\infty)\|u^{k/2}\|_2^2
  \end{array}$$
  by H\"older's inequality.
 On the other hand, we have
$$ \begin{array}{rcl} |((u_k\cdot\na)u_{k/2,k},u^{k})|
 &\leq &  \|u_k\|_\infty\|\na u_{k/2,k}\|_2\|u^k\|_2\ek
  &\leq &  k\|u_k\|_\infty\|u^{k/2}\|_2^2.
 \end{array}
 $$
 Therefore, we get from \eq{2.18} that
 \begin{equation}
 \label{E2.19}
 \begin{array}{l}
 \frac{d}{2dt}\|u^k\|_2^2+\nu\|\na u^k\|_2^2
      \leq (\|\na u_{k/2}\|_\infty+\|\na u_k\|_\infty+k\|u_k\|_\infty)\|u^{k/2}\|_2^2
      \quad\text{in }(0,T).
      \end{array}
 \end{equation}

Let
 \begin{equation}
 \label{E4.3}
j_0\equiv j_0(k)=\lceil\log_2 k\rceil+1.
 \end{equation}
Since
 $$\begin{array}{rcl} \|\widehat{\Da_ju_k}\|_1
 &\leq & \|\widehat{\Da_ju_k}\|_2\cdot |\{2^{j-1}\leq |\xi|\leq
 2^j\}|^{1/2}\ek
  &\leq &  c_02^{3j/2}\|\Da_ju_k\|_2\ek
   &\leq &  c_02^{3j/2}\|\Da_ju\|_2,\;\forall j\in\Z,
 \end{array}$$
 holds with a generic constant $c_0>0$, we get by
  the definition of the space $X_1$
 that
 $$\begin{array}{rl}
 \|u_k\|_\infty&\leq (2\pi)^{-3/2}\|\widehat{u_k}\|_1\leq (2\pi)^{-3/2}
  \di\sum_{j=-\infty}^{j_0}\|\widehat{\Da_ju_k}\|_1\ek
  &\leq c_0(2\pi)^{-3/2}\di\sum_{j=-\infty}^{j_0} 2^{j}a(j)
     (a^{-1}(j)2^{j/2}\big\|\Da_ju\|_2)\ek
 &\leq
 \di c_0(2\pi)^{-3/2}\big(\sum_{j=-\infty}^{j_0} 2^{j}a(j)\big)
   \cdot \sup_{j\leq j_0}\{a^{-1}(j)2^{j/2}\big\| \Da_ju \big\|_2\}\ek
 &\leq c_0(2\pi)^{-3/2}(2+  b(j_0(k)) 2^{j_0+1})  \big\|u\big\|_{X_1}\ek
   &\leq 8c_0(2\pi)^{-3/2}b(j_0(k))k\|u\|_{X_1}
  \end{array}$$
  due to $k\geq 1, b(j)\geq 1/2$ for all $j\in\N$, see \eq{2.2} and \eq{3.2}.
In the same way, we have
 $$\begin{array}{rl}
 \|\na u_k\|_\infty &\leq
 \di c_0(2\pi)^{-3/2}\big(\sum_{j=-\infty}^{j_0} 2^{j}a(j)\big)
   \cdot \sup_{j\leq j_0}\{a^{-1}(j)2^{j/2}\big\| \Da_j\na u \big\|_2\}\ek
 &\leq \di c_0(2\pi)^{-3/2}(2+  b(j_0(k)) 2^{j_0+1})\cdot 2^{j_0}
   \sup_{j\leq j_0}\{a^{-1}(j)2^{j/2}\big\| \Da_ju \big\|_2\}\ek
   &\leq 8c_0(2\pi)^{-3/2}b(j_0(k))k^2\|u\|_{X_1}
  \end{array}$$
and
 $$\begin{array}{rl}
 \|\na u_{k/2}\|_\infty &\leq
 \di c_0(2\pi)^{-3/2}\big(\sum_{j=-\infty}^{j_0-1} 2^{j}a(j)\big)
   \cdot \sup_{j\leq j_0-1}\{a^{-1}(j)2^{j/2}\big\| \Da_j\na u \big\|_2\}\ek
&\leq
 \di c_0(2\pi)^{-3/2}\big(\sum_{j=-\infty}^{j_0} 2^{j}a(j)\big)
   \cdot 2^{j_0-1}\sup_{j\leq j_0-1}\{a^{-1}(j)2^{j/2}\big\| \Da_j u \big\|_2\}\ek
   &\leq 4c_0(2\pi)^{-3/2}b(j_0(k))k^2\|u\|_{X_1}.
  \end{array}$$

 Thus,  we have by \eq{2.19} that
 \begin{equation}
 \label{E2.15}
\begin{array}{l}
\di  \frac{d}{2dt}\|u^{k}\|_2^2+\nu \|\na u^k\|_2^2
   \leq  C_1b(j_0(k))k^2\|u\|_{X_1}\|u^{k/2}\|_2^2\ek
\qquad\leq  4C_1b(j_0(k))\|u\|_{X_1}\|\na u^{k/2}\|_2^2,  \,\forall
k\geq 1, \quad\text{in }(0,T),
\end{array}
 \end{equation}
with a generic constant $C_1>0$ independent of $\nu$.

{\it Step 2. Higher order norm estimate:}

Let us add up \eq{2.15} over all $k\in \N$. Then we have
 \begin{equation}
 \label{E2.5}
  \begin{array}{l}
 \di\frac{d}{2dt}\sum_{k\in\N}\|u^{k}\|_2^2+\nu \sum_{k\in\N}
 \|\na u^{k}\|_2^2  \leq  4C_1\|u(t)\|_{X_1}  \sum_{k\in\N} b(j_0(k))\|\na u^{k/2}\|_2^2\;\text{in }(0,T).
   \end{array}
 \end{equation}
Note that
 \begin{equation}
 \label{E2.25}
\begin{array}{l} \di \sum_{k\in \N}\|u^{k}\|_2^2
=\di\sum_{n\in\N}n\|u_{n,n+1}\|_2^2,\quad
 \di \sum_{k\in\N}\|\na u^{k}\|_2^2 =\di\sum_{n\in\N}n\|\na u_{n,n+1}\|_2^2.
 \end{array}
 \end{equation}
Moreover, with the notation $2\N:=\{2n: n\in\N\}$ we have
 \begin{equation}
 \label{E2.26}
\begin{array}{l}
  \di\sum_{k\in\N} b(j_0(k))\|\na u^{k/2}\|_2^2\\[3ex]
   \di = \ha\|\na u^{1/2}\|_2^2+
   \sum_{k\in 2\N} \big(b(j_0(k))\|\na u^{k/2}\|_2^2+b(j_0(k+1))\|\na u^{(k+1)/2}\|_2^2\big)
\\[3ex]
   \di \leq  \|\na u^{1/2}\|_2^2+
   \sum_{k\in 2\N} \big(b(j_0(k))+b(j_0(k+1))\big)\|\na u^{k/2}\|_2^2
\\[3ex]
   \di =  \|\na u^{1/2}\|_2^2+
   \sum_{k\in \N} (b(j_0(2k))+b(j_0(2k+1)))\| \na u^{k}\|_2^2
\\[3ex]
%   \di \leq  \| u^{1/2}\|_2^2+9
%   \sum_{k\in \N} \big(b(j_0(2k))+b(j_0(2k+1))\big)k^2\| u^{k}\|_2^2
%   \\[3ex]
\di =\|\na u^{1/2}\|_2^2+\sum_{k\in\N}
\big(b(j_0(2k))+b(j_0(2k+1))\big)\sum_{n\geq k}\|\na u_{n,n+1}\|_2^2\\[3ex]
\di =\|\na u^{1/2}\|_2^2+\sum_{n\in\N} \|\na u_{n,n+1}\|_2^2
\sum_{k=1}^n \big(b(j_0(2k))+b(j_0(2k+1))\big).
   \end{array}
  \end{equation}
Note that $j_0(2k)=\lceil \log_2k\rceil +2$, see \eq{4.3}, and
 $$\begin{array}{l}
\di\sum_{k=1}^n b(j_0(2k))\\[3ex]
 \quad= \di\sum_{k\leq n,j_0(2k)\notin S(\lceil\log_2n\rceil+2)} b(j_0(2k))
 +\sum_{k\leq n,j_0(2k)\in S(\lceil\log_2n\rceil+2)} b(j_0(2k)).
 \end{array}$$
where, by Lemma \ref{L4.1},
 $$\begin{array}{l}
\di\sum_{k\leq n,j_0(2k)\notin S(\lceil\log_2n\rceil+2)} b(j_0(2k))
\leq 2n.
 \end{array}$$
Moreover,  we have by \eq{4.1} and Lemma \ref{L4.2} that
$$\begin{array}{l}
\di\sum_{k\leq n,j_0(2k)\in S(\lceil\log_2n\rceil+2)} b(j_0(2k))\\[3ex]
 \qquad \qquad \di\leq  |S(\lceil\log_2n\rceil+2)|\cdot
           \max_{k\leq n, j_0(2k)\in S(\lceil\log_2n\rceil+2)} b(j_0(2k))\\[3ex]
   \qquad \qquad  \di\leq  |S(\lceil\log_2n\rceil+2)|\cdot
   \log_2(\lceil\log_2n\rceil+2)\\[3ex]
\qquad \qquad  \di\leq
\frac{3\log_2\log_2(\lceil\log_2n\rceil+2)+5}{2}\cdot
\log_2\log_2(\lceil\log_2n\rceil+2)\cdot\log_2(\lceil\log_2n\rceil+2)\\[3ex]
\qquad \qquad  \di \leq n,\quad \forall n\geq n_0,
 \end{array}$$
 for some generic number $n_0\in\N$.
Therefore, for all $n\geq n_0$ we have
$$\begin{array}{l}
\di\sum_{k=1}^n b(j_0(2k)) \leq 2n+n\leq 3n.
 \end{array}$$
In the same way, we have
$$\begin{array}{l}
\di\sum_{k=1}^n b(j_0(2k+1)) \leq 3n,\quad \forall n\in\N, n\geq
n_0,
 \end{array}$$
with a generic number $n_0\in\N$  possibly larger than the former $n_0$.

Thus we get by \eq{2.26} that
 \begin{equation}
 \label{E2.26n}
\begin{array}{l}
  \di\sum_{k\in\N} b(j_0(k))\| \na u^{k/2}\|_2^2\leq \|\na u^{1/2}\|_2^2\\[3ex]
\qquad\di +\sum_{n=1}^{n_0-1} \|\na u_{n,n+1}\|_2^2 \sum_{k=1}^n
\big(b(j_0(2k))+b(j_0(2k+1))\big)\\[3ex]
\qquad  +\di\sum_{n\geq n_0} \|\na u_{n,n+1}\|_2^2 \sum_{k=1}^n \big(b(j_0(2k))+b(j_0(2k+1))\big)\\[3ex]
 \quad\di \leq \|\na u^{1/2}\|_2^2+(c(n_0)-1)\|\na u_{n_0}\|_2^2+6\sum_{n\geq n_0}n\|\na u_{n,n+1}\|_2^2
\\[3ex]
 \quad\di \leq c(n_0)\|\na u\|_2^2+6\sum_{n\geq n_0} n\|\na u_{n,n+1}\|_2^2,
   \end{array}
  \end{equation}
  where $c(n_0)$ is a generic constant depending only on $n_0$.

Therefore, we get from \eq{2.5}, \eq{2.25} and \eq{2.26n} that
 \begin{equation}
  \label{E4.13}
  \begin{array}{l}
 \di\frac{d}{2dt}\di\sum_{n\in\N}n\|u_{n,n+1}\|_2^2+\nu \sum_{n\in\N}
 n\|\na u_{n, n+1}(t)\|_2^2  \\[3ex]
 \quad\leq 4C_1\|u(t)\|_{X_1}\big(c(n_0)\|\na u(t)\|_2^2
   +6\di \sum_{n\geq n_0} n\|\na u_{n,n+1}(t)\|_2^2\big)\ek
\quad \leq C_2\|u(t)\|_{X_1}\big(\|\na u(t)\|_2^2
   +\di \sum_{n\in\N} n\|\na u_{n,n+1}(t)\|_2^2\big)\quad \text{in
   }(0,T),
   \end{array}
 \end{equation}
where  $C_2>0$ is a generic constant.
%Thus, if
% \begin{equation}
% \label{E2.20}
%18C_1m < \frac{\nu}{2},
% \end{equation}
%then
% \begin{equation}
%  \label{E4.14}
% \begin{array}{l}
% \di\frac{d}{dt}\di\sum_{n\in\N}n\|u_{n,n+1}\|_2^2+\nu \sum_{n\in\N}
% n\|\na u_{n, n+1}\|_2^2  \leq 2C_1m \tilde{c}(n_0)\|\na u\|_2^2.
%   \end{array}
%   \end{equation}

By integrating \eq{4.13} from $0$ to $t\in (0,T)$  we have
 \begin{equation}
 \label{E4.12}
 \begin{array}{l}
\di\ha\sum_{n\in\N}n\|u_{n,n+1}(t)\|_2^2+\nu \int_0^t \sum_{n\in\N}
 n\|\na u_{n, n+1}(\tau)\|_2^2\,d\tau
\leq \ha\di\sum_{n\in\N}n\|[u_0]_{n,n+1}\|_2^2\\[2ex]
  \quad +C_2\di\int_0^t \|u(\tau)\|_{X_1}\big(\|\na u(\tau)\|_2^2
   +\di \sum_{n\in\N} n\|\na u_{n,n+1}(\tau)\|_2^2\big)\,d\tau,
   \quad \forall t\in (0,T).
  \end{array}
  \end{equation}
Since
 $$n\|u_{n,n+1}\|_2^2\leq \int_{n\leq |\xi|\leq
n+1}|\xi||\hat{u}(\xi)|^2d\xi\leq (n+1)\|u_{n,n+1}\|_2^2,$$ we have
$$\begin{array}{l}
\di \sum_{n\in\N}n\|u_{n,n+1}\|_2^2\leq
\|u^1\|_{\dot{H}^{1/2}}^2\leq
\sum_{n\in\N}(n+1)\|u_{n,n+1}\|_2^2=\sum_{n\in\N}n\|u_{n,n+1}\|_2^2+\|u^1(t)\|_2^2,\ek
 \di\sum_{n\in\N}n\|\na
u_{n,n+1}\|_2^2\leq \|u^1\|_{\dot{H}^{3/2}}^2\leq
\sum_{n\in\N}n\|\na u_{n,n+1}\|_2^2+\|\na u^1(t)\|_2^2,
  \end{array}$$
and the estimate \eq{4.12} implies that
 \begin{equation}
 \label{E4.14}
 \begin{array}{l}
\di\ha\|u^1(t)\|_{\dot{H}^{1/2}}^2+\nu \int_0^{t}
 \|u^1(\tau)\|_{\dot{H}^{3/2}}^2\,d\tau
\leq \ha\|u^1(t)\|_{2}^2+\nu \int_0^{t}
 \|\na u^1(\tau)\|_{2}^2\,d\tau+\ha\|u_0\|_{\dot{H}^{1/2}}^2\\[2ex]
  \quad +\di C_2\int_0^t \|u(\tau)\|_{X_1}\big(\|\na u(\tau)\|_2^2
   + \|u^1(\tau)\|_{\dot{H}^{3/2}}^2\big)\,d\tau,
   \quad \forall t\in (0,T).
  \end{array}
  \end{equation}
Since
$\|u\|_{\dot{H}^s}^2=\|u_1\|_{\dot{H}^s}^2+\|u^1\|_{\dot{H}^s}^2,
s\geq 0,$ and
 \begin{equation}
 \label{E4.2}
\ha\|u_1(t)\|_{\dot{H}^{1/2}}^2+\nu \int_0^{t}
 \|u_1(\tau)\|_{\dot{H}^{3/2}}^2\,d\tau
\leq \ha\|u_1(t)\|_{2}^2+\nu \int_0^{t}
 \|\na u_1(\tau)\|_{2}^2\,d\tau,
  \end{equation}
we get by summing \eq{4.14} and \eq{4.2} that
 \begin{equation}
 \label{E4.14n}
 \begin{array}{l}
\di\ha\|u(t)\|_{\dot{H}^{1/2}}^2+\nu \int_0^{t}
 \|u(\tau)\|_{\dot{H}^{3/2}}^2\,d\tau
\leq \ha\|u(t)\|_{2}^2+\nu \int_0^{t}
 \|\na u(\tau)\|_{2}^2\,d\tau+\ha\|u_0\|_{\dot{H}^{1/2}}^2\\[2ex]
  \quad +\di C_2\int_0^t \|u(\tau)\|_{X_1}\big(\|\na u(\tau)\|_2^2
   + \|u(\tau)\|_{\dot{H}^{3/2}}^2\big)\,d\tau\\[2ex]
   \leq \di \ha\|u_0\|_{2}^2+\ha\|u_0\|_{\dot{H}^{1/2}}^2
 +C_2\int_0^t \|u(\tau)\|_{X_1}\big(\|\na u(\tau)\|_2^2
   + \|u(\tau)\|_{\dot{H}^{3/2}}^2\big)\,d\tau,\forall t\in (0,T).
  \end{array}
  \end{equation}

Now, put $\ve=\frac{\nu}{2C_2}$ and fix $t\in (0,T)$. By Lemma
\ref{L4.2n} (ii) in view of $u\in C([0,t],
  \dot{H}^{1/2}(\R^3))$ there exists
$l_0(t)>0$ such that
 \begin{equation}
 \label{E4.16}
\begin{array}{l}
  \di\|\la u(\tau, \la\cdot)\|_{X_1}\leq \ve,\quad  \forall \tau\in [0,t],\;\forall
\la=2^{2^{2^l}} (l\geq l_0(t)).
 \end{array}
 \end{equation}

  Now, put $\la=2^{2^{2^l}}$ with fixed $l\geq l_0(t)$  and observe
  the re-scaled function
 \begin{equation}
 \label{E4.19}
 (u)_\la(\tau,x):=\la u(\la^2 \tau, \la x), \,\tau\in
 (0,\la^{-2}t].
  \end{equation}
Then, \eq{4.16} implies that
 \begin{equation}
 \label{E4.17}
\begin{array}{l}
 \di \|(u)_\la(\tau)\|_{X_1}\leq \ve,\quad  \forall \tau\in [0,\la^{-2}t].
 \end{array}
 \end{equation}

Since $(u)_\la$ is the Leray-Hopf weak solution to \eq{1.1} with
initial value $(u_0)_\la$ with the first blow-up epoch $\la^{-2}T$,
we have by \eq{4.14n} that
 \begin{equation}
 \label{E4.14nn}
 \begin{array}{l}
\di\ha\|(u)_\la(\la^{-2}t)\|_{\dot{H}^{1/2}}^2+\nu
\int_0^{\la^{-2}t}
 \|(u)_\la(\tau)\|_{\dot{H}^{3/2}}^2\,d\tau
\leq \ha\|(u_0)_\la\|_2^2+\ha\|(u_0)_\la\|_{\dot{H}^{1/2}}^2\\[2ex]
  \quad +\di C_2\int_0^{\la^{-2}t} \|(u)_\la(\tau)\|_{X_1}\big(\|\na (u)_\la(\tau)\|_2^2
   +\|(u)_\la(\tau)\|_{\dot{H}^{3/2}}^2\big)\,d\tau.
  \end{array}
  \end{equation}
By \eq{4.17}, \eq{4.14nn} we have
 \begin{equation}
 \label{E4.18}
 \begin{array}{l}
\di\ha\|(u)_\la(\la^{-2}t)\|_{\dot{H}^{1/2}}^2+\frac{\nu}{2}
\int_0^{\la^{-2}t}
 \|(u)_\la(\tau)\|_{\dot{H}^{3/2}}^2\,d\tau
\\[2ex]
  \qquad \di\leq \ha\|(u_0)_\la\|_2^2+\ha\|(u_0)_\la\|_{\dot{H}^{1/2}}^2+\di\nu \int_0^{\la^{-2}t} \|\na
  (u)_\la(\tau)\|_2^2\,d\tau
  \\[2ex]
  \qquad \di\leq
  \|(u_0)_\la\|_2^2+\ha\|(u_0)_\la\|_{\dot{H}^{1/2}}^2,
  \end{array}
  \end{equation}
  where we used energy inequality for $(u)_\la$, i.e.,
$$\di\ha\|(u)_\la(\la^{-2}t)\|_2^2+\nu \int_0^{\la^{-2}t} \|\na
  (u)_\la(\tau)\|_2^2\,d\tau\leq \ha\|(u_0)_\la\|_2^2.$$

Thanks to the scaling properties of the norms, that is,
$$\begin{array}{l}
\|(u)_\la(\la^{-2}t)\|_{\dot{H}^{1/2}}=\|u(t)\|_{\dot{H}^{1/2}},\;
 \|(u_0)_\la\|_{\dot{H}^{1/2}}=\|u_0\|_{\dot{H}^{1/2}},\ek
 \di\int_0^{\la^{-2}t}
 \|(u)_\la(\tau)\|_{\dot{H}^{3/2}}^2\,d\tau=\int_0^{t}
 \|u(\tau)\|_{\dot{H}^{3/2}}^2\,d\tau,\;\|(u_0)_\la\|_2=\la^{-1/2}\|u_0\|_2,
  \end{array}$$
it follows from \eq{4.18} that
\begin{equation}
 \label{E4.22}
 \begin{array}{rcl}
\di\ha\|u(t)\|_{\dot{H}^{1/2}}^2+\frac{\nu}{2}\int_0^{t}
 \|u(\tau)\|_{\dot{H}^{3/2}}^2\,d\tau
  &\leq &  \la^{-1}\|u_0\|_{2}^2+\ha\|u_0\|_{\dot{H}^{1/2}}^2\\[2ex]
  &\leq &  \ha\|u_0\|_2^2+\ha\|u_0\|_{\dot{H}^{1/2}}^2\\[2ex]
   &\leq &  \ha\|u_0\|_{H^{1/2}}^2
  \end{array}
  \end{equation}
  due to $\la>2$.
 %Note that the right-hand side of \eq{4.22} is independent of $t\in (0,T)$.
 Since \eq{4.22} holds for any fixed $t\in (0,T)$,  we conclude that $T=\infty$.
 Moreover, \eq{4.22} and the energy inequality \eq{EI} yield \eq{1.2} and
$$u\in L^\infty(0,\infty; H^{1/2}(\R^3))\cap  L^2(0,\infty;
H^{3/2}(\R^3)).
$$
 The proof is complete.

\hfill\qed

\bigskip
%\noindent {\bf Fund:} This work was supported by no fund organization.

\end{document}